\documentclass[12pt]{amsart}
\usepackage{amssymb}
\usepackage[T1]{fontenc}
\usepackage{amsmath,stackrel}
\usepackage{indentfirst}
\usepackage{xspace}
\usepackage{color}
\usepackage{fancyhdr} 
\usepackage{fancybox}
\usepackage[brazil,english]{babel}
\usepackage[latin1]{inputenc}
\usepackage{graphicx}
\usepackage[all]{xy}
\usepackage{hyperref}

\theoremstyle{plain}
\newtheorem{teo}{Theorem}
\newtheorem{lema}[teo]{Lemma}
\newtheorem{prop}[teo]{Proposition}
\newtheorem{cor}[teo]{Corollary}

\theoremstyle{definition}
\newtheorem{defi}{Definition}[section]
\newtheorem{exe}{Example}

\theoremstyle{remark}
\newtheorem*{obs}{Remark}
\numberwithin{teo}{section}

\newcommand{\C}{\mathbb{C}} 

\newenvironment{prova}
{{\em{\noindent \textbf{Proof:} }}
}
{\hfill $\square$}

\DeclareMathOperator{\rank}{rank}
\DeclareMathOperator{\im}{im}
\DeclareMathOperator{\codim}{codim}

\begin{document}

\title[Euler obstruction of essentially isolated determinantal singularities]{Euler obstruction of essentially isolated determinantal singularities}
\author{Nancy Carolina Chachapoyas Siesqu\'en}

{\small\address{Universidade Federal de Itajub\'{a}, Instituto de Matem\'{a}tica e Computa\c{c}\~{a}o.
Avenida BPS, 1303, Pinheirinho, CEP: 37500903 - Itajub\'{a}- MG, Brasil}}
\email{nancy@unifei.edu.br}

\date{}
\begin{abstract}
  We study the Euler obstruction of essentially isolated determinantal singularities (EIDS). The EIDS were defined  by W. Ebeling and S. Gusein-Zade \cite{Gusein}, as a generalization of isolated singularity. We obtain some formulas to calculate the Euler obs\-truction for the determinantal varieties with singular set an ICIS.  

\end{abstract}

\maketitle

\thispagestyle{plain}
\pagestyle{fancy}
\fancyhead{}
\fancyhead[CO]{\bfseries \small{Euler obstruction of essentially isolated determinantal singularities}}
\fancyhead[LE]{Nancy Carolina Chachapoyas Siesqu\'en }
\renewcommand{\headrulewidth}{0pt}


\section{Introduction}
\label{sec:introduction}

\thispagestyle{empty}
In this work, we study the Euler obstruction of essentially isolated determinantal singularities (EIDS).  The EIDS have been defined by W.~{E}beling and S.~M. Guse{\u\i}n-Zade in \cite{Gusein}. A generic determinantal variety $M_{m,n}^t$ is the subset of the space of $m \times n$ matrices, given by the matrices of rank less than 
$t$, where $0 \leqslant t\leqslant \min\{m,n\}$.
 A variety $X\subset \C^N$ is  determinantal if $X$ is the pre-image of $M_{m,n}^t$ by a holomorphic function $F:\C^N\to M_{m,n}$ with the condition that the codimension of $X$ in $\C^N$ is the same as the codimension of $M_{m,n}^t$ in $M_{m,n}$.

Determinantal varieties have isolated singularities if $N\leqslant (m-t+2)(n-t+2)$  and they admit smoothing if $N<(m-t+2)(n-t+2)$, see \cite{Wahl}. Seve\-ral recent works investigate determinantal varieties with isolated singularities. The Milnor number of a determinantal surface was defined in \cite{Pike, Bruna, MIRIAMSP} while the vanishing Euler characteristic of a determinantal variety was defined in \cite{Pike, Bruna}. Other recent results on isolated determinantal varieties related to this paper  appear in particular in \cite{Anne-Matthias, Gaffney-Antoni, Gaffney-Ruas}. Section of EIDS are studied in \cite{BCR}.

In this work we study  the Euler obstruction of EIDS, specially for EIDS admitting a stratification with at most 3 strata.
The main result is a formula to compute the Euler obstruction in terms of the singular vanishing Euler characteristic.

\section{Essentially isolated determinantal singularity }

We denote by $M_{m,n}$ the set of matrices $m\times n$ with complex entries.
\begin{defi} For all  $t$, $ 1\leqslant t\leqslant \min\{m,n\}$, let $M_{m,n}^{t}$ be  the subset of $M_{m,n}$ whose elements are matrices of rank less than $t$:
$$M_{m,n}^{t}=\{A\in M_{m,n} \vert  \rank(A)< t \}.$$
This set is a singular variety of codimension $ (m-t+1)(n-t+1) $ in $M_{m,n}$, called generic determinantal variety.
\end{defi}
The singular set of $M_{m,n}^{t}$ is $M_{m,n}^{t-1}$.
The partition of $M_{m,n}^{t}$ defined by 
$$M_{m,n}^{t}=\cup_{i=1,...,t}( M_{m,n}^{i} \backslash M_{m,n}^{i-1})$$ is a Whitney stratification \cite{Harris}.

Let $F: \C^{N}\rightarrow M_{m,n}$ be a map defined by $F(x)=(f_{ij}(x))$, whose entries are complex analytic functions defined on an open domain  $U\subset \C^N$.

\begin{defi} The map $F$ determines a determinantal variety of type $(m,n,t)$ in $U\subset \C^N$ as the analytic variety $X=F^{-1}(M_{m,n}^{t})$, such that $\codim X=\codim M_{m,n}^t=(m-t+1)(n-t+1)$.
\end{defi}

A generic map $F$ intersects transversally the strata $M_{m,n}^i \backslash M_{m,n}^{i-1}$ of the variety $M_{m,n}^t$. The following definition was introduced in \cite{Gusein}.
\begin{defi} A point $x\in X= F^{-1}(M_{m,n}^{t})$ is called essentially nonsingular if at this point the map $F$ is transversal to the corres\-ponding stratum of the variety $M_{m,n}^{t}$ (that is, to $\{M_{m,n}^{i} \backslash M_{m,n}^{i-1}\}$, where $i=\rank F(x)+1$ ).
\end{defi}

%
%
%

\begin{defi}\label{EIDS} A germ $(X,0)\subset (\C^{N},0) $ of a determinantal variety of type $(m,n,t)$ has an essentially isolated singular point at the origin (or is an essentially isolated determinantal singularity, EIDS) if it has only essentially nonsingular points in a punctured neighborhood of the origin in $X$.
\end{defi}

An EIDS $X\subset \C^N$ has isolated singularity  if and only if $N\leqslant (m-t+2)(n-t+2)$. An EIDS with isolated singularity will be called isolated determinantal singularity, denoted by IDS \cite{Bruna, MIRIAMSP}.

We want to consider deformations of an EIDS that are themselves determinantal varieties of the same type.

\begin{defi} An essential smoothing $\widetilde{X}_s$ of the EIDS $(X,0)$ is a subvariety lying in a neighborhood of the origin in $\C^N$ and defined by $\widetilde{X}_s = \widetilde{F}^{-1}_s (M_{m,n}^t)$ where 
 $\widetilde{F}:U\times \C \to M_{m,n}$ is a perturbation of the germ $F$, with $\widetilde{F}_{s}(x)=\widetilde{F}(x,s)$, $\widetilde{F}_0(x)=F(x)$ such that $\widetilde{F}_s:U\to M_{m,n}$ is transversal to all strata $M_{m,n}^i \backslash M_{m,n}^{i-1}$.
\end{defi}
An essential smoothing is in general not smooth (when $ N\geq (m-t+2)(n-t+2)$) as we see in the following theorem.

\begin{teo}\cite{Wahl} \label{smoothing} Let $(X,0)\subset (\C^N,0)$ be the germ of a determinantal variety with isolated singularity at the origin. Then, $X$ has a smoothing if and only if  $N<(m-t+2)(n-t+2).$
\end{teo}

When the essential smoothing $\widetilde{X}_s$ is singular, its singular set is $\widetilde{F}_s^{-1}(M^{t-1}_{m,n})$. Since $\widetilde{F}$ is transversal to the strata of the Whitney stratification $M_{m,n}^t$, the partition  $\widetilde{X}_s=\cup_{1\leq i \leq t} \widetilde{F}_s^{-1}(M^{i}_{m,n}\backslash M^{i-1}_{m,n})$ is a Whit\-ney stratification of $\widetilde{X}_s$.


\begin{exe} Let $X=F^{-1}(M_{2,3}^2)$ be the $4$-dimensional variety in $\C^6$ where
 $$\begin{array}{cccl}
                F :& \C^{6}      & \rightarrow & M_{2,3} \\
                   & (x,y,z,w,v, u)   & \mapsto     &\left(
                                                 \begin{array}{ccc}
                                                   x & y& v \\
                                                   z & w   & x+u^2 \\
                                                 \end{array}
                                               \right).

              \end{array}
$$
The following matrix defines an essential smoothing $\widetilde{X}_s=\widetilde{F}_s^{-1}(M_{2,3}^2)$ of $X$
$$\begin{array}{cccl}
                \widetilde{F} :& \C^{6}  \times \C   & \rightarrow & M_{2,3} \\
                   & (x,y,z,w,v, u,s)   & \mapsto     &\left(
    \begin{array}{ccc}
      x+s & y & v \\
      z & w & x+u^2 \\
    \end{array}
  \right).
              \end{array}$$
              
 In this case, $\widetilde X_s$ is singular.

\end{exe}

\section{L\^e-Greuel formula type for IDS with smoo\-thing}

In this section we review some results about isolated singularities following \cite{MIRIAMSP}.

Let $(X,0)\subset (\C^N,0)$ be the germ of a $d$- dimensional variety with isolated singularity at the origin. Suppose that $X$ has a smoothing. Then, there exists a flat family $\pi: \widetilde{X}\subset U\times \C\to \C$ such that the fiber $X_s=\pi^{-1}(s)$ is smooth for all $s\neq 0$ and $X_0=X$.

Let $p:(X,0)\to \C$ be a complex analytic function defined in $X$ with isolated singularity at the origin. Let us consider a function
$$\begin{array}{cccc}
      \widetilde{p}:&   \C^N\times \C   &   \to &  \C \\
      & (x,s)  & \mapsto & \widetilde{p}(x,s),
\end{array}$$
such that $ \widetilde{p}(x,0)=p(x)$ and for all $s\neq 0$, $ \widetilde{p}(\cdot,s)=p_s$ is a Morse function on $X_s.$

Thus we have the following diagram

\begin{displaymath}
\xymatrix{ X_s \ar[d]^{p_s} \ar@{^{(}->}[r]  &  \C^N\times \C \ar[d]^{(\pi,\widetilde{p})}\\
                    \C\times \{s\}    \ar@{^{(}->}[r]     & \C\times \C}
\end{displaymath}

\begin{prop} \cite[Proposition 4.1]{MIRIAMSP}  \label{prop:Miriam} Let $X$ be a d-dimensional variety with isolated singularity at the origin admitting smoothing and $p_s:X_s \to \C$, $p_s=\widetilde{p}(\cdot,s)$ as above.
Then,

\begin{itemize}
\item [(a)] If $s\neq0$, $X_s\simeq p_s^{-1}(0)\dot{\cup}\{ \text{cells of dimension d } \},$

\item [(b)] $\chi(X_s)=\chi(p_s^{-1}(0))+(-1)^d n_{0},$
\end{itemize}
where  $n_{0}$ is the number of critical points of $p_s$ and $\chi(X_s)$ denotes the Euler characteristic of $X_s$.
\end{prop}

Let us recall how the invariant $n_{0}$ is related to the polar multiplicity of $X$, $m_d(X)$ (\cite{MIRIAMSP}, see also \cite{Gaffney1}).
\begin{defi}(The $d$-Polar multiplicity) Let ${X}$, $\widetilde{X}$, $p$ and $\widetilde{p}$ as above.
Let $P_d(X,\pi,{p})=\overline{\Sigma(\pi,\widetilde{p})|}_{\widetilde{X}_{reg}}$ be the relative polar variety of $X$ related to $\pi$ and $p$. We define $m_d(X,\pi,{p})=m_0(P_d(X,\pi,{p}))$.
\end{defi}

In general $m_d(X,\pi,{p})$ depends on the choices of $\widetilde{X}$ and $\widetilde{p}$. When the variety $X$ has a 
unique smoothing $\widetilde{X}$, then 
$m_d(X,\pi,{p})$ depends only on $X$ and ${p}$. If ${p}$ is a generic linear embedding, $m_d(X,p)$  is an invariant of the EIDS $X$, denoted by $m_d(X)$.

\begin{prop} \cite{MIRIAMSP} Under the conditions of Proposition \ref{prop:Miriam}, $n_{0}=m_d(X)$.
\end{prop}

\begin{teo}  \cite{greuellivro} Let $X_s$ be a smoothing of a normal isolated singularity, then $b_1(X_s)=0.$
\end{teo}

Let $X$ be a determinantal variety of type $(m,n,t)$ in $\C^N$ with $N< (m-t+2)(n-t+2)$. Then $X$ has isolated singularity and admits smoothing. 
\begin{defi} \cite{MIRIAMSP}\label{Defi:Milnor} Let $X$ be a determinantal surface in $\C^N$, with isolated singularity at the origin. The Milnor number of $X$, denoted by $\mu(X) $, is defined as the second Betti number of the generic fiber $X_s$,
$$\mu(X)=b_2(X_s).$$
\end{defi}

The following result appears in \cite{Pike,Bruna, MIRIAMSP}, for determinantal surfaces $X\subset \C^4$, but it also holds for any surface with isolated singularity in $\C^N$ admitting smoothing.

\begin{prop}\label{LeGreuel} Let $(X,0)\subset (\C^N,0)$ be the germ of a determinantal surface in $\C^N$ with isolated singularity at the origin admitting smoo\-thing. Let $p:(\C^N,0)\to (\C,0)$  be a linear function whose restriction to $X$ has an isolated singularity at the origin. Then one has the L\^e-Greuel formula
\begin{equation}
\label{formulaGreuel}
\mu(X)+\mu(X\cap p^{-1}(0))=m_2(X,p).
\end{equation}
\end{prop}

When $d=\dim X> 2,$ the Betti numbers $b_i(X),$ $2 \leqslant i < d$ are not ne\-ces\-sarily zero (see \cite{Anne-Matthias}). In \cite{Pike, Bruna} the authors define the vanishing Euler charac\-teristic of varieties admitting smoothing.

\begin{defi} \cite{Bruna} Let $(X,0)\subset (\C^N,0)$ be an IDS such that $N < (m-t+2)(n-t+2)$. The vanishing Euler characteristic is defined by \[\nu(X)= (-1)^d(\chi(X_s)-1),\]
where $X_s$ is a smoothing of  $X$ and $\chi(X_s)$ is the Euler characteristic of $X_s$ .
\end{defi}

 \begin{teo} \label{Greuel} \cite{Bruna} Let $(X,0)\subset (\C^N,0)$ be an IDS such that $N<(m-t+2)(n-t+2)$ and let $p:(\C^N,0)\to (\C,0)$ be a linear projection whose restriction to $X$ has isolated singularity at the origin. Then,
\begin{equation}
\label{}
\nu(X)+\nu(X\cap p^{-1}(0))=m_d(X,p).
\end{equation}
When $p$ is a generic linear projection, then $m_d(X,p)=m_d(X)$.
\end{teo}

\begin{obs} When $d=2$, then $\nu(X)=\mu(X)$.
\end{obs}

\begin{exe} \cite{MP} Let $X=F^{-1}(M_{2,3}^2)\subset \C^4$ be the variety defined by:
$$\begin{array}{cccc}
      F:&  \C^4& \to &M_{2,3}\\
      &(x,y,z,w)   & \mapsto & \begin{pmatrix}
    z  & y+w  & x \\
    w  &  x & y
\end{pmatrix}.
\end{array}$$

The following matrix defines an essential smoothing $\widetilde{X}_s=\widetilde{F}_s^{-1}(M_{2,3}^2)$ of $X$
$$\begin{array}{cccl}
                \widetilde{F} :& \C^{4}  \times \C   & \rightarrow & M_{2,3} \\
                   & (x,y,z,w,s)   & \mapsto     &\left(
    \begin{array}{ccc}
      z & y+w & x+s \\
      w & x & y \\
    \end{array}
  \right).
              \end{array}$$
  In this case $\widetilde X$ is an smoothing of $X$.
Consider $p:\C^4\to \C$ given by $p(x,y,z,w)=w$, then it follows that $m_2(X)=3$ and $\mu(X\cap p^{-1}(0))=2$, then $\mu(X)=1.$
\end{exe}

\section{Euler obstruction}
Let $(X,0)\subset (\C^N,0)$ be the germ of a complex analytic variety. Let $Gr(d,N)$ be the Grassmanian of complex $d$-planes in $\C^N$. On the regular part $X_{reg}$ of $X$ the Gauss map $\phi:X_{reg}\to \C^N\times Gr(d,N)$ is well defined by $\phi(x)=(x,T_x(X_{reg}))$, where $T_x(X_{reg})$ denotes the tangent space of $X_{reg}$ at the point $x$.
\begin{defi} \label{DefiNash} The Nash transformation $\widetilde{X}$ of $X$ is the closure of the image $\im \phi$ in $\C^N\times Gr(d,N)$. That is,
 $$\widetilde{X}=\overline{\{(x, W)/ x\in X_{reg},  W=T_{x}X_{reg}\}}\subset X\times Gr(d,m).$$
 
It is a complex analytic space endowed with an analytic projection map $\nu: \widetilde{X} \to X$ which is a biholomorphism away from $\nu^{-1}(Sing(X))$.
 \end{defi}
Denote by $E$  the tautological bundle over $Gr(d, N)$.  We denote still by $E$ the corresponding trivial extension bundle over $ \C^N \times  Gr(d,N)$ with projection map $\pi$. 
We define the Nash bundle  $\widetilde{T}$ on $\widetilde{X}$, as the restriction of $E$ to $\widetilde{X}$.
\begin{displaymath}
\xymatrix{  \widetilde{T}  \ar@{^{(}->}[r]  \ar[d]^{\pi}&E \ar[d]\\
                  \widetilde{X} \ar@{^{(}->}[r] \ar[d]^{\nu}& \C^N \times  Gr(d,N) \ar[d]  \\
                    X \ar@{^{(}->}[r]        & M.}
\end{displaymath}

Let  $X$ be a small representative of the germ $(X,0)$ and $\{V_i\}$ a Whitney stratification of $X$. Suppose that  the point$\{0\}$ is a stratum. Let $\cup T V_{i}$ be the union of the tangent bundles of all the strata. This union can be seen as a subset of the tangent bundle $T\C^N|_X.$ We follow the references \cite{Brasselet, Brasselet2, Schwartz}.

\begin{defi} A stratified vector field $v$ on $X$ is a continuous section of  $T\C^N|_X$ such that, if $x\in V_i\cap X$, then $v(x) \in T_x(V_i)$. 
\end{defi}

\begin{defi} A radial vector field $v$ in a neighbourhood of $\{0\}$ in $X$ is a stratified vector field so that there is $\epsilon_0$ such that for every $0< \epsilon \leqslant \epsilon_0$, $v(x)$ is pointing outwards the ball   $B_{\epsilon}$ over the boundary $S_{\epsilon}= \partial B_{\epsilon}$
\end{defi}

\begin{lema} \cite{Schwartz} Let $X$ be a stratified non-zero vector field on $Z\subset X$. Then $v$ can be lifted as a section $\widetilde v$ of $\widetilde T$ over $\nu^{-1}(Z)$.

\end{lema}

The definition of Euler obstruction was given by  R. MacPherson \cite{MACPHERSON}. Here we give an equivalent definition given by J.-P. Brasselet and M.-H. Schwartz.
\begin{defi} \cite{Schwartz} Let $v$ be a radial vector field on $X\cap \partial B_{\epsilon}$ and $\tilde{v}$ the lifting of $v$ on $\nu^{-1}(X\cap \partial B_{\epsilon}).$ The vector $\tilde{v}$ defines a cocycle of obstruction $Obs(\tilde{v})$, that measures the obstruction to extending $\tilde{v}$ as a non zero section of $\tilde{T}$ over $\nu^{-1}(X\cap B_{\epsilon}):$
$$Obs(\tilde{v})\in Z^{2N}(\nu^{-1}(X\cap B_{\epsilon}), \nu^{-1}(V\cap \partial B_{\epsilon}))$$
The local Euler  obstruction (or Euler  obstruction ) denoted by $Eu_X(0)$ is the evaluation of the cocycle $Obs(\tilde{v})$ over the fundamental class of $X$, $[\nu^{-1}(X\cap B_{\epsilon}), \nu^{-1}(V\cap \partial B_{\epsilon})]$, that is,

$$Eu_0(X)=\langle Obs(\tilde{v}), [\nu^{-1}(X\cap B_{\epsilon}), \nu^{-1}(V\cap \partial B_{\epsilon})]\rangle$$
\end{defi}

\begin{obs} The following properties of the Euler obstruction hold
\begin{enumerate}
\item $Eu_x(X)$ is constant along the strata of a Whitney stratification of $X.$
\item $Eu_x(X)=1$, if $x\in X_{reg}$.
\item Assume that $X$ is a local embedding in $\C^N$ and $W$ a smooth variety $\C^N$ that intersects the Whitney stratification of  $X$ transversely. Then $Eu_x(X\cap W)=Eu_x(X)$ for all $x\in X\cap W$.
\item $Eu_{(x,y)}(X\times Y)= Eu_x(X)\times Eu_y(Y)$ for $x\in X$ and $y\in Y.$
\end{enumerate}
\end{obs}

The following theorem gives a Lefschetz type formula for the local Euler obstruction \cite[theorem 3.1]{Brasselet2} that is important tool to compute the Euler obstruction.

\begin{teo}\label{formula} \cite{Brasselet2} Let $(X,0)$ be a germ of an equidimensional complex analytic space in $\C^N$. Let $V_i,$ $i=1,\dots,l$ be the strata of the Whitney stratification of a small representative $X$ of $(X,0)$. Then there is an open dense Zariski subset $\Omega$ in the space of complex linear forms on $\C^N$, such that for every $l\in \Omega$, there is $\epsilon_0,$ such that for any $\epsilon,\epsilon_0>\epsilon>0$ and  $t_0\neq 0$ sufficiently small, such that the following formula for the Euler obstruction holds 
$$Eu_0(X)=\sum_{1}^{l}\chi(V_i\cap B_{\epsilon}\cap l^{-1}(t_0))Eu_{V_i}(X),$$
where $Eu_{V_i}(X)$ is the value of the Euler obstruction of $X$ at any point of $V_i,$ $i=1,\dots,l$.
\end{teo}
To compute the Euler obstruction of determinantal varieties, the following formula given in \cite{Gusein}  is an useful tool. For determinantal varieties  $X=F^{-1}(M_{m,n}^t)$ outside the origin,  a normal section to the stratum $V_i=X_i\backslash X_{i-1} $ (where $X_i=F^{-1}(M_{m,n}^i)$) is isomorphic to $(M_{m-i+1,n-i+1}^{t-i+1},0).$
\begin{teo}\label{Matrix} \cite{Gusein} Let $l:M_{m,n}\to \C$ be a generic linear form and let $L_{m,n}^t=M_{m,n}^t\cap l^{-1}(1)$. Then, for $t\leqslant m \leqslant n$.
$$\overline{\chi}(L_{m,n}^t)=(-1)^t\begin{pmatrix}
      m-1  \\
      t-1  
\end{pmatrix},$$
where $\overline{\chi}(L_{m,n}^t)=\chi(L_{m,n}^t)-1.$
\end{teo}

\section{ Singular Vanishing Euler Characteristic of Determinantal Varieties}

In the case that $V$ is not a complete intersection, the authors in \cite{Pike} define the singular vanishing Euler characteristic as follows:
\begin{defi} \label{deficarac} \cite{Pike} Let $V\subset \C^p$ be a complex analytic variety (not necessarily a complete intersection). Suppose that $F:\C^N \to \C^p$ is transverse to  the strata of $V$ outside the point $\{0\}$ and consider $F_s:\C^N \to \C^p$ a 1-parameter deformation of $F$ such that $F_s$ is transversal to $V$ for all $s\neq 0$ small. We say that $F_s$ is a "stabilization''  of $F$. Write $X=F^{-1}(V)$. Then the singular vanishing Euler characteristic of $X$ is given by:
$$\tilde{\chi}(X)= \tilde{\chi}(F_s^{-1}(V))=\chi(F_s^{-1}(V))-1, \hspace{0.2cm}. $$
\end{defi}

\begin{obs} \cite{Pike} If $V$ is a $k$-dimensional complete intersection and $X=F^{-1}(V)$ then $$\tilde{\chi}(X)=(-1)^{N-k}\mu(X).$$
\end{obs}

In the work  \cite{Pike}, the authors obtain several formulas to calculate the singular vanishing Euler characteristic of determinantal varieties defined by $2\times 3$ matrizes. Proposition \ref{Pike} is a Corollary of Theorem 8.1 in \cite{Pike}.

If $F:\C^N\to M_{2,3}$ is a generic germ of corank $1$, then $$F(x,y)=\sum_{i=1}^{5} x_iw_i+g(y)w_0,$$
where $\{w_1,w_2,w_3,w_4,w_5\}$ is a basis for $W=dF(0)(\C^N)$, $w_0\notin W$ and $g(y)$ define a isolated singularity on $\C^{N-5}.$

\begin{prop}  \cite[Corollary 11.5]{Pike} \label{Pike}   If $F$  is a generic germ of corank $1$, $N\geqslant 6$, 
$$\tilde{\chi}(X)=(-1)^{N-1}\mu(g).$$
If $g$ is quasi-homogeneous, $\mu(g)$ is $\tau(X).$
\end{prop}

\section{Euler obstruction of determinantal varieties in $\C^N$}

Let $X=F^{-1}(M_{m,n}^t)$ be an EIDS, defined by $F:\C^N\to M_{m,n}$. If $N\leqslant (m-t+3)(n-t+3)$ then the singular set  $\Sigma X=F^{-1}(M_{m,n}^{t-1})$  is a determinantal variety with isolated singularity. Hence, the variety $X$ admits at most $3$  strata $\{V_0,V_1 ,V_2\}$, where $V_0=\{0\}$, $V_1=\Sigma X\backslash \{0\}$, $V_2=X_{reg}$. These are the determinantal varieties we consider in this section.

For the calculus of the Euler obstruction we use the formula of the Theorem \ref{formula}, then we have
\begin{eqnarray} \nonumber
Eu_0(X) & = & \chi(V_0\cap l^{-1}(r)\cap B_{\epsilon})Eu_{V_0}(X)+ \chi(V_1\cap l^{-1}(r) \cap B_{\epsilon})Eu_{V_1}(X) \\ \nonumber
 & & + \chi(V_2\cap l^{-1}(r)\cap B_{\epsilon})Eu_{V_2}(X).  \nonumber
\end{eqnarray}
As $V_0\cap l^{-1}(r)\cap B_{\epsilon}=\emptyset$ then $\chi(V_0\cap l^{-1}(r)\cap B_{\epsilon})=0$. We also have $Eu_{V_2}(X)=1$, then
\begin{equation}
\label{eq:Eu}
Eu_0(X) =  \chi(V_1\cap l^{-1}(r)\cap B_{\epsilon})(Eu_{V_1}(X)-1)+\chi(X\cap l^{-1}(r)\cap B_{\epsilon}) 
\end{equation}

Note that $ \chi(V_1\cap l^{-1}(r)\cap B_{\epsilon})=  \chi(\Sigma X \cap l^{-1}(r)\cap B_{\epsilon})$. Then, given $a\in V_1$ there is an open set $U_a$ containing  $a$ such that $U_a\cong B_{\delta}\times c(L_{V_1})$, where the dimension of $B_{\delta}$ is $\dim V_1=N-(m-t+2)(n-t+2)$ and $c(L_{V_1})$ is the cone over the complex link of $V_1$ in $X$, $L_{V_1}=X\cap \mathcal{N} \cap p^{-1}(s)$, with $\mathcal{N}$ transversal to $V_1$ at $a$, $\codim \mathcal{N}=\dim V_1$, then 
\begin{eqnarray} \nonumber
Eu_{V_1}(X) & = & Eu_a (B_{\delta})\times Eu_a(c(L_{V_1})) \\ \nonumber
 & = &  Eu_a(c(L_{V_1}))=\chi(L_{V_1}) \\ \nonumber 
  & = & \chi(X\cap \mathcal{N}\cap p^{-1}(s)). \nonumber
\end{eqnarray}
The equality $Eu_a(c(L_{V_1}))=\chi(L_{V_1}) $ is given  by applying again Theo\-rem \ref{formula} to $c(L_{V_1})$, observing that $c(L_{V_1})$ has a isolated singularity and $c(L_{V_1})\cap p^{-1}(s)$ is isomorphic to $L_{V_1}$.
Here $X\cap \mathcal{N}$ is an essential smoothing of a determinantal variety of type $(m,n,t)$ in $\C^r$, where $r=\codim V_1=(m-t+2)(n-t+2)$, therefore with isolated singularity.

Substituting  $Eu_{V_1}(X)$ in \eqref{eq:Eu}, we have
\begin{equation}
\label{ }
Eu_0(X)  =    \chi(\Sigma X \cap l^{-1}(r)\cap B_{\epsilon})(\chi(L_{V_1})-1)+\chi(X\cap l^{-1}(r)\cap B_{\epsilon}).
\end{equation}
This formula can be expressed in terms of the singular vanishing Euler characteristic of the Definition \ref{deficarac}.
\begin{equation}
\label{geral}
Eu_0(X)  =    (\tilde{\chi}(\Sigma X \cap l^{-1}(0)\cap B_{\epsilon})+1)(\chi(L_{V_1})-1)+\tilde{\chi}(X\cap l^{-1}(0)\cap B_{\epsilon})+1.  \nonumber
\end{equation}

\begin{prop} \label{ICISOBS} Let $X=F^{-1}(M_{m,n}^t)$be an EIDS, defined by $F:\C^N\to M_{m,n}$. If $N\leqslant(n-t+3)(m-t+3)$ and $\Sigma X$ is an ICIS, then
\begin{eqnarray} 
 Eu_0(X) & = & ((-1)^{\dim (\Sigma X\cap l^{-1}(0) )}\mu(\Sigma X\cap l^{-1}(0))+1)(\chi(L_{V_1})-1) \\
 &  & +\tilde{\chi}(X\cap l^{-1}(0)\cap B_{\epsilon})+1 \nonumber
\end{eqnarray}
where $l:\C^N \to \C$ is a generic linear projection, $L_{V_1}$ is the complex link of the stratum $V_1$ em $X$ and $B_{\epsilon}$ is the ball of radius ${\epsilon}$ in $\C^N$.
\end{prop}

\begin{prova}

If $\Sigma X$ is an ICIS, then $$\chi(\Sigma X\cap l^{-1}(r))=1+(-1)^{\dim (\Sigma X\cap l^{-1}(0))}\mu(\Sigma X\cap l^{-1}(0)).$$ Substituting this formula in \eqref{geral} we have the result.
\end{prova}

The following result shows the existence of determinantal varieties whose singular set is an ICIS. 

\begin{prop} Let $X=F^{-1}(M_{2,n}^2)\subset \C^N$ be an EIDS defined by the function $F:\C^N\to M_{2,n}$, $n\geqslant 2$ then the singular set of  $X$, $\Sigma X\subset \C^N$ is an ICIS.
\end{prop}

\begin{prova} As $X=F^{-1}(M_{2,n}^2)$ is an EIDS then $\Sigma X=F^{-1}(M_{2,n}^1)$. As $M_{2,n}^1=\{0\}$, then $\Sigma X=F^{-1}(0)$ and $\codim \Sigma X= (2-2+2)(n-2+2)=2\cdot n$, therefore $\Sigma X$ is an ICIS.

Let $x\in \Sigma X$, $x\neq 0$ then $\rank F(x)=0$, by the definition of EIDS, we have that $F$ is transversal to the corresponding stratum, in this case $F$ is transversal to $\{0\}$ at the origin. Then $dF$  is an embedding in $x\neq 0$, so $\Sigma X$ is an ICIS.
\end{prova}


\subsection{ Euler obstruction of IDS, case $N<(n-t+2)(m-t+2)$}
In this section we will treat the varieties IDS,  defined by the function $F:\C^N\to M_{m,n}$, such that $N<(n-t+2)(m-t+2)$. These varieties admit smoothing and the Euler obstruction was studied in the work \cite{Bruna}, using the vanishing Euler characteristic and the multiplicity of the variety.

\begin{teo} \cite{Bruna} \label{Bruna} Let $X=F^{-1}(M_{m,n}^t)$ be the determinantal variety defined by  $F:\C^N\to M_{m,n}$ with $N<(n-t+2)(m-t+2)$, then
$$Eu_0(X)=1+(-1)^d\nu(X,0)+(-1)^{d+1}m_d(X,0)$$
\end{teo}
\begin{exe} Let $X=F^{-1}(M_{2,3}^2)\subset \C^4$ be a determinantal variety with isolated singularity, defined by $F$.
$$\begin{array}{cccl}
      F:& \C^4  &\to& M_{2,3} \\
      &(x,y,z,w)   &\mapsto & \begin{pmatrix}
      x &  y & z \\
      y & z & w
\end{pmatrix}
\end{array}$$
By \cite{Bruna} we have that $m_2(X)=3$, $\nu(X)=1$, then $Eu_0(X)=-1.$
\end{exe}

\subsection{Euler obstruction of an EIDS, case $N=6$}

The next result is for the case $N=6$. In this case the determinant variety does not admit smoothing. The formula of the Theorem \ref{Bruna}, given above does not hold.
\begin{teo} Let $X=F^{-1}(M_{2,3}^2)$ defined by the $F:\C^N\to M_{2,3}$ with $N= 6$, then
$$Eu_0(X)=b_2(X\cap l^{-1}(r))-b_3(X\cap l^{-1}(r))+1.$$
\end{teo}
\begin{prova} The formula comes from the fact that $$\chi(X\cap l^{-1}(r))-1=\tilde{\chi}(X\cap l^{-1}(0))=b_2(X\cap l^{-1}(r))-b_3(X\cap l^{-1}(r)).$$
\end{prova}

In the Tables \ref{Table1} and \ref{Table2}, we calculate the Euler obstruction of $4$- dimensional determinantal varieties in $\C^6$, classified in \cite{NEUMER}. We use the Table $5$ of article \cite{Pike} to obtain the results. The invariant $\tau$ in Tables \ref{Table1} and \ref{Table2} is the Tjurina number.
\begin{table}
  \centering 
  \begin{tabular}{ccccc}
\hline
\footnotesize{Type} & $\begin{array}{c}\text{\footnotesize{ Form of }}\\ \text{\footnotesize{the matrix}} \end{array}$& \footnotesize{Conditions} & $\tau$ & $Eu_0(X)$  \\
\hline \hline 
& & & &\\
$ \Omega_1$&  $  \begin{pmatrix}
   x   & y  & v \\
   z & w  & u
\end{pmatrix}$ & & 0 & 2 \\
  & & & &\\ 
 $ \Omega_k$&  $\begin{pmatrix}
   x   & y  & v \\
   z & w  & x+u^k
\end{pmatrix}$  & $k\geqslant 2$ &$k-1$ & 2\\
& & & &\\
   $A_k^\dag$&  $\begin{pmatrix}
      x & y & z    \\
    w  & v & u^2+x^{k+1}+y^2  
\end{pmatrix}$& $k\geqslant 1$ & $k-2$ & 1 \\
& & & &\\
$D_k^\dag$& $\begin{pmatrix}
     x & y & z    \\
     w & v & u^2+xy^2+x^{k-1}  
\end{pmatrix}$ & $k\geqslant 4$ & $k+2$& -1\\
& & & &\\
$E_6^\dag$& $\begin{pmatrix}
    x  & y & z   \\
     w &  v & u^2+y^4
\end{pmatrix}$ & & 8 & 0\\
& & & &\\
$E_7^\dag$& $\begin{pmatrix}
     x & y & z    \\
      w & v & u^2+x^3+xy^3  
\end{pmatrix}$ & & 9 & 0\\
& & & &\\
$E_8^\dag$&$\begin{pmatrix}
      x & y & z    \\
     w & v & u^2+x^3+y^5  
\end{pmatrix}$ & &$10$ & 0\\
& & & &\\
&$\begin{pmatrix}
     x & y & z   \\
     w & v & ux+y^k+u^l  
\end{pmatrix}$ & $k\geqslant 2$, $l\geqslant 3$& $k+l-1$ &$3-k$\\
& & & &\\
&$\begin{pmatrix}
  x    & y & z    \\
  w    & v & x^2+y^2+u^3  
\end{pmatrix} $& & $6$ & $1$\\
& & & &\\
$F_{q,r}^\dag$& $\begin{pmatrix}
     w & y & x    \\
     z & w+vu & y+v^q+u^r  
\end{pmatrix}$ & $q, r\geqslant 2$ & $q+r$& $2$ \\

& & & &\\
$G_5^\dag$& $\begin{pmatrix}
     w & y & x    \\
  z    & w+v^2 & y+u^3 
\end{pmatrix}$& & $7$ & $2$\\

& & & &\\
$G_7^\dag$& $\begin{pmatrix}
     w & y & x    \\
   z   & w+v^2 & y+u^4 
\end{pmatrix}$ & & $10$ & $2$\\
& & & &\\
\hline
\end{tabular}
  
  \caption{ Euler obstruction, $X\subset \C^6$}\label{Table1}
\end{table}

\begin{table}
  \centering 
  
   \begin{tabular}{llccc}
\hline
\footnotesize{Type} & $\begin{array}{c}\text{\footnotesize{ Form of}}\\ \text{\footnotesize{the matrix}} \end{array}$& \footnotesize{Conditions} & $\tau$ & $Eu_0(X)$  \\
\hline \hline 
& & & &\\

$H_{q+3}^\dag$& $\begin{pmatrix}
    w  &  y & x  \\
    z & w+v^2+u^q & y+vu^2  
\end{pmatrix}$ & $q\geqslant 3$ &  $q+5$& 2\\

& & & &\\
$ I_{2q-1}^\dag$&  $  \begin{pmatrix}
   w  & y  & v \\
   z & w+v^2+u^3  & y+u^q
\end{pmatrix}$ & $q\geqslant4$ & $2q+1$ & 2 \\
  & & & &\\ 
$ I_{2r+2}^\dag$&  $\begin{pmatrix}
   w   & y  & x \\
   z & w+v^2+u^3  & y+vu^r
\end{pmatrix}$  & $r\geqslant 3$ &$2r+4$ & 2\\
& & & &\\
  & $\begin{pmatrix}
     w & y & x    \\
     z & w+v^2 & u^2+yv  
\end{pmatrix}$ &   & $6$& 1\\
& & & &\\
& $\begin{pmatrix}
    w  & y & x   \\
     z &  w+uv & u^2+yv+v^k
\end{pmatrix}$ &$k\geqslant3$ & $k+4$ & $1$ \\
& & & &\\
& $\begin{pmatrix}
     w & y & x    \\
      z & w+v^k & u^2+yv+v^3  
\end{pmatrix}$ &$k\geqslant3  $ & $2k+2$ & 1\\
& & & &\\
&$\begin{pmatrix}
      w & y & x    \\
     z & w+uv^k & u^2+yv+v^3  
\end{pmatrix}$ & $k\geqslant2$&$2k+5$ & 1\\
& & & &\\
  & $\begin{pmatrix}
     w & y & x    \\
  z  & w+v^3& u^2+yv 
\end{pmatrix}$& & 9 & $1$\\
& & & &\\
 & $\begin{pmatrix}
     w & y & x    \\
   z   & w+v^k & u^2+y^2=v^3 
\end{pmatrix}$ & $k\geqslant3$& $2k+3$ & $1$\\
& & & &\\
$H_{q+3}^\dag$& $\begin{pmatrix}
    w  &  y & x  \\
    z & w+uv^k & u^2+y^2+v^3  
\end{pmatrix}$ & $k\geqslant 2$ &  $2k+6$& 1\\
& & & &\\
\hline
\end{tabular}
  \caption{Euler obstruction, $X\subset \C^6$.}\label{Table2}
\end{table}

\subsection{ Euler obstruction of EIDS, case $N\geqslant 7$}

\begin{teo} \label{N7} Let $X\subset\C^N$ be an EIDS, defined by the function $F:\C^N\to M_{2,3}$, with $N\geqslant 7$. Then\begin{equation}
\label{N71}
Eu_0(X)=(-1)^{N-7}\mu(\Sigma X\cap l^{-1}(0))+\tilde{\chi}(X\cap l^{-1}(0))+2.
\end{equation}
\end{teo}
\begin{prova} In this case, $X$ has three strata, $\{V_0,V_1,V_2\}$, $V_0=\{0\}$, $V_1=\Sigma \backslash \{0\}$, $V_2=X_{reg}$ and $V_1=F^{-1}(M_{2,3}^1)\backslash F^{-1}(M_{2,3}^0)$.  We have that $L_{V_1}=X\cap \mathcal{N}\cap p^{-1}(s)$, where $X\cap \mathcal{N}\cong M_{2,3}^2$. So $\chi(L_{V_1})-1=1$ by Theorem \ref{Matrix}. Substituting this  value in the Proposition \ref{ICISOBS} we have the result.
\end{prova}
%
%
%

\begin{cor} \label{2} With the hypotheses of the Theorem \ref{N7}, if $F$ has corank $1$, then $$Eu_0(X)=2.$$
\end{cor}
\begin{prova}
As $F$ have corank $1$, then $$F(x,y)=\sum_{i=1}^{5} x_iw_i+g(y)w_0,$$
where $\{w_1,w_2,w_3,w_4,w_5\}$ is a base for $W=d_0F(\C^N)$, $w_0\notin W$ and $g(y_1,y_2, \dots, y_{N-5})$ define a isolated singularity on $\C^{N-5}.$ Suppose wi\-thout loss of generality that $l(x,y)=y_1$ and consider $\tilde{g}(y_2, \dots, y_{N-5})=g(0,y_2,\dots, y_{N-5})$

 Using the formula of the Proposition \ref{Pike} for $X\cap l^{-1}(0)\subset \C^{N-1}$,we have that $\tilde{\chi}(X\cap l^{-1}(0))=(-1)^{N-2}\mu(\tilde{g})$. 
 
 Also, we have that $$\Sigma X\cap l^{-1}(0)=\{(0,\dots, 0, y_2,\dots, y_{N-5}) | \tilde{g}(y_2,\dots, y_{N-5})=0\}$$ and hence $\mu(\Sigma X\cap l^{-1}(0))= \mu(\tilde{g}).$
Substituting in the formula \eqref{N71}, we have that $Eu_0(X)=2.$
\end{prova}

\begin{exe} Let $X$ be the variety defined by the function $F$.
$$\begin{array}{cccc}
      F:&\C^8 & \to  & M_{2,3} \\
      &(x,y)  &\mapsto & \begin{pmatrix}
     x_1 & x_2 &x_3   \\
      x_4& x_5 & x_1+y_1^2-y_2^2+y_3^2  
\end{pmatrix}
\end{array}$$
where $x=(x_1,x_2,x_3,x_4,x_5)$ e $y=(y_1,y_2,y_3)$. The variety $X$ is an EIDS of dimension $6$ in $\C^8$. The singular set has dimension $2$ and is the set $$ \Sigma X=\{(0,0,0,0,0,y_1,y_2,y_3)| y_1^2-y_2^2+y_3^2=0  \}.$$
This variety has three strata $V_0=0,$ $V_1=\Sigma X-\{0\}$ and $V_2=X_{reg}$.
 
 As the corank of $F$ is 1, then by the Corollary \ref{2}, we have 
  $$Eu_0(X)= 2.$$
 \end{exe}

%
 \section*{Acknowledgments}
 The author would like to thank J.-P. Brasselet and  M. A. S Ruas for their careful reading and suggestions to improve this article.
 This work was partially supported by FAPESP grants no. 2010/ 09736-1, 2011/20082-6 and CNPq grant no. 164353/2014-3.
 
\nocite{Nancy}

\end{document}